\documentclass{amsart}
\usepackage[utf8]{inputenc}
\usepackage{amsmath,amssymb,color,mathrsfs,enumerate}
\usepackage{float}
\usepackage{graphicx}
\usepackage[pagebackref]{hyperref}

\usepackage{xcolor}
\usepackage{comment}
\usepackage{mathtools}
\usepackage[new]{old-arrows}

\numberwithin{equation}{section}
\theoremstyle{plain}
\newtheorem{theorem}{Theorem}[section]
\newtheorem{lemma}[theorem]{Lemma}
\newtheorem{proposition}[theorem]{Proposition}
\newtheorem{corollary}[theorem]{Corollary}

\theoremstyle{definition}
\newtheorem{definition}[theorem]{Definition}

\newtheorem{example}[theorem]{Example}

\newtheorem{notation}[theorem]{Notation}

\newtheorem{para}[theorem]{}

\newtheorem{setup}[theorem]{Setup}

\theoremstyle{remark}
\newtheorem{remark}[theorem]{Remark}

\newcommand{\N}{\mathbb{N}}

\newcommand{\Z}{\mathbb{Z}}

\newcommand{\mycomment}[1]{}

\DeclareMathOperator{\ann}{ann}
\DeclareMathOperator{\Ass}{Ass}

\DeclareMathOperator{\eend}{end}

\DeclareMathOperator{\indeg}{indeg}
\DeclareMathOperator{\ini}{in}
\DeclareMathOperator{\init}{in_<}

\DeclareMathOperator{\Mon}{Mon}

\DeclareMathOperator{\rank}{rank}
\DeclareMathOperator{\reg}{reg}

\renewcommand{\ge}{\geqslant}
\renewcommand{\le}{\leqslant}

\newcommand{\bz}{\mathbb{Z}}

\newcommand{\fa}{\mathfrak{a}}
\newcommand{\fm}{\mathfrak{m}}

\newcommand{\fp}{\mathfrak{p}}
\newcommand{\fq}{\mathfrak{q}}

\newcommand{\fr}{\mathfrak{r}}

\definecolor{MyDarkGreen}{cmyk}{0.7,0,1,0}

\title[Vasconcelos invariants for products and powers of ideals]{Asymptotic behaviour of Vasconcelos invariants for products and powers of graded ideals}

\author[L.~Fiorindo]{Luca Fiorindo} 
\address{Dipartimento di Matematica, Dipartimento di Eccellenza 2023-2027, Università di Genova, Via Dodecaneso 35, 16146 Genova, Italy}
\email{luca.fiorindo@edu.unige.it}
\urladdr{\url{https://orcid.org/0000-0002-6435-0128}}

\author[D.~Ghosh]{Dipankar Ghosh}
\address{Department of Mathematics, Indian Institute of Technology Kharagpur, West Bengal - 721302, India}
\email{dipankar@maths.iitkgp.ac.in, dipug23@gmail.com}
\urladdr{\url{https://orcid.org/0000-0002-3773-4003}}



\begin{document}

\begin{abstract}
Let $R$ be a commutative Noetherian $\N$-graded ring. Let $N\subseteq M$ be finitely generated $\Z$-graded $R$-modules. Let $I_1,\ldots,I_r$ be nonzero proper homogeneous ideals of $R$. Denote ${\bf I}^{\underline{n}}:=I_1^{n_1}\cdots I_r^{n_r}$ for $\underline{n}=(n_1,\dots,n_r)\in\N^r$. In this paper, we prove that the (local) Vasconcelos invariant of ${\bf I}^{\underline{n}}M/{\bf I}^{\underline{n}}N$ is eventually the minimum of finitely many linear functions in $\underline{n}$. The same holds for $M/{\bf I}^{\underline{n}}N$ under certain conditions. Some specific examples are provided, where these functions are not eventually linear in $\underline{n}$. However, when $R$ is a polynomial ring over a field, we show that the global Vasconcelos invariants of $R/{\bf I}^{\underline{n}}$ and ${\bf I}^{\underline{n}}/{\bf I}^{\underline{n}+\underline{1}}$ are, in fact, asymptotically linear in $\underline{n}$ with the leading coefficients given by the initial degrees of $I_1,\ldots,I_r$. The last result is surprising: It differs from the Castelnuovo-Mumford regularity, which is not always linear even over polynomial rings, as shown by Bruns-Conca.
\end{abstract}
\maketitle

\section{Introduction}
The Vasconcelos invariant, also known in the literature as $v$-number, is a recent invariant known both in the communities of Commutative Algebra and Coding Theory. It takes its name from the mathematician Wolmer V. Vasconcelos (1937-2021). It was first introduced in \cite[Sec.~4.1]{CSTPV} for a homogeneous ideal $I$ of a polynomial ring $R$ over a field: For an associated prime $\fp\in\Ass(R/I)$, the local $v$-number is the non-negative integer
$$v_\fp(I):=\inf\{n\ge 0: \text{there exists } f\in R_n\text{ such that }\fp=(I:_R f)\}.$$
Moreover, $v(I) :=\inf\{ v_\fp(I) : \fp\in\Ass(R/I)\}$ is called the $v$-number of $I$. This invariant was first used to express the regularity index of the minimum distance function of projective Reed-Muller-type codes. In the following years, mathematicians have studied this invariant in different areas spanning between Commutative Algebra and Combinatorics, discovering connections with other different invariants. 
For instance, a combinatorial interpretation is shown in \cite[Thm.~A]{JS} for the Vasconcelos invariant of a binomial edge ideal using the connected domination number. A similar result is also obtained in \cite[Thm.~3.5]{JV} for square-free monomial ideals. 
In \cite[p.~16]{CSTPV}, the authors connect the Vasconcelos invariant to the degree of finite projective varieties. Various (in)equalities between Castelnuovo-Mumford regularity and $v$-number are established in \cite[Thm.~4.10]{CSTPV}, \cite[Thm.~3.13]{JV}, \cite[p.~905]{SS}, \cite[Thm.~3.8]{Sa}, \cite[Thm.~4.19]{BM}, \cite[Prop.~2.2 and Rmk.~2.3]{FG} and \cite[Thm.~A]{NS}.

A different approach is given in \cite{FS} by Ficarra-Sgroi considering the following problem: Since $\Ass(R/I^n)$ stabilizes to a set $\Ass^\infty(R/I)$ for $n$ big enough, due to a result of Brodmann \cite{brod}, it is possible to study the behaviour of $v(I^n)$ as a function of $n$. In \cite{FS}, it is proved that the functions $v_\fp(I^n)$ for $\fp\in \Ass^\infty(R/I)$, and $v(I^n)$ are eventually linear; moreover, the leading coefficient of $v(I^n)$ is the initial degree of $I$. Instead, independently in \cite{Conca23}, Conca proved the same result in a more general frame when $R$ is a Noetherian standard graded domain, and he showed that the leading coefficient of $v_\fp(I^n)$ lies in the degrees in which the ideal $I$ is generated. Recently, in \cite{FG}, the authors extended the notion of Vasconcelos invariant to a finitely generated graded module $M$ over a Noetherian graded ring $R$. For $\fp\in\Ass_R(M)$, the local $v$-number of $M$ is given by
\[
v_\fp(M):=\inf\{n\in\Z: \text{there exists } x\in M_n\text{ such that }\fp=(0:_R x)\},
\]
while the quantity $v(M):=\inf\{v_\fp(M):\fp\in\Ass(M)\}$ is called the $v$-number of $M$, see \cite[Defn.~1.2]{FG}. By convention, $v(0) =\infty$. For a homogeneous ideal $I$ of $R$, when $(0:_M I)=0$, it is shown in \cite[Thm.~1.9]{FG} that $v(M/I^{n}M)$ is eventually linear, where the leading coefficient is explicitly described. This considerably strengthens the results of Conca and Ficarra-Sgroi. See \cite[Thm.~3.8]{GP25} for a more general result. The main aim of the present article is to consider the products and powers of several homogeneous ideals, that is ${\bf I}^{\underline{n}}:=I_1^{n_1}\cdots I_r^{n_r}$, and to study the asymptotic behaviour of the corresponding Vasconcelos invariant. The motivation for our results came from the asymptotic behaviour of the Castelnuovo-Mumford regularity of ${\bf I}^{\underline{n}}M$ as a function in $\underline{n}$. When $R$ is a standard graded algebra over a field, in \cite[Cor.~4.4]{Gh16}, it is shown that $\reg({\bf I}^{\underline{n}}M)$ is bounded above by a linear function in $\underline{n}$. Later, Bruns-Conca in \cite[Thm.~2.2]{BC17} proved that asymptotically $\reg({\bf I}^{\underline{n}}M)$ is, in fact, the maximum of finitely many linear functions in $\underline{n}$.

To better present the results of this paper, we fix the following notations.

\begin{setup}\label{setup}
    Let $R$ be a commutative Noetherian $\N$-graded ring. Let $M$ be a finitely generated $\mathbb{Z}$-graded $R$-module. For each $1\le i\le r$, suppose $I_i$ is a homogeneous ideal of $R$ generated in degrees $d_{i,j}$ for $1 \le j \le a_i$.
    Let $N$ be a graded submodule of $M$ (e.g., $N=\mathfrak{a}M$ for some homogeneous ideal $\mathfrak{a}$ of $R$). Set ${\bf I} := I_1\cdots I_r$. For $\underline{n}=(n_1,\dots,n_r)\in\N^r$, denote ${\bf I}^{\underline{n}}:=I_1^{n_1}\cdots I_r^{n_r}$.
\end{setup}

\begin{para}\label{para:notations-r-tuple}
    In this paper, the additive group $\Z^r$, of $r$-tuples $\underline{n}=(n_1,\dots,n_r)$ of integers with componentwise addition, is endowed with the componentwise order, that is $\underline{n}\ge \underline{m}$ if $n_i\ge m_i$ for all $i=1,\dots,r$. By $\underline{0}$ and $\underline{1}$, we denote the $r$-tuples $(0,\dots,0)$ and $(1,\dots,1)$ respectively. Let $\underline{e}_j$ for $1\le j\le r$ denote the standard basis of $\Z^r$ as a free $\Z$-module. For $\underline{m}, \underline{n}\in\Z^r$, set $\underline{m}\cdot \underline{n} := m_1 n_1 + \cdots + m_r n_r$, which is the usual dot product of $\underline{m}$ and $\underline{n}$. For a property (P), by writing ``(P) holds true for all $\underline{n}\gg \underline{0}$", we mean ``there exists $\underline{m}\in\N^r$ such that (P) holds true for all $\underline{n}\ge \underline{m}$".
\end{para}

\begin{para}\label{para:notations}
    For an ideal $I$ of $R$, we use the notations:
    \begin{center}
    $(N:_MI) := \{x\in M:Ix\subseteq N\}, \quad \Gamma_I(M) := \bigcup_{n \ge 1}\big(0:_M I^n\big)$,
    \end{center}
    $\ann_M(I) := (0:_M I)$, and $\ann_R(M) := \{ r \in R : rM = 0\}$. The initial degree of $M$ is defined to be $ \indeg(M) := \inf\{ n\in\Z : M_n \not=0 \} $. By convention, $\indeg(0) := +\infty$.
\end{para}


\begin{notation}\label{notation}
    With Setup~\ref{setup}, by \cite[Cor.~1.2]{Hay06}, the sets $\Ass_R(M/{\bf I}^{\underline{n}}N)$ and $\Ass_R({\bf I}^{\underline{n}}M/{\bf I}^{\underline{n}}N)$ stabilize (possibly to two different sets) for all $\underline{n}\gg\underline{0}$. So we denote $\mathcal{A}_N^M({\bf I}) := \Ass_R(M/{\bf I}^{\underline{n}}N)$ and $\mathcal{B}_N^M({\bf I}) := \Ass_R({\bf I}^{\underline{n}}M/{\bf I}^{\underline{n}}N)$ for all $\underline{n}\gg\underline{0}$.
\end{notation}


Our main results are summarized in the following two theorems.

\begin{theorem}\label{thm::asymbehaviour}
    With {\rm Setup~\ref{setup}} and {\rm Notation \ref{notation}}, the following statements hold.
    \begin{enumerate}[\rm (1)]
        \item For each $\fp\in \mathcal{B}_{N}^M({\bf I})$, there exist $\underline{w}_1,\dots,\underline{w}_s\in\N^r$ and $c_1,\dots,c_s\in\Z$ such that
        $$v_\fp({\bf I}^{\underline{n}}M/{\bf I}^{\underline{n}} N)=\min\{\underline{w}_k\cdot\underline{n}+c_k:\, 1\le k\le s\}$$
        for $\underline{n}\gg\underline{0}$. Moreover, if $\underline{w}_{k} = (w_{k1},w_{k2},\ldots,w_{kr})$, then $w_{ki}\in\{d_{i,1},\dots,d_{i,a_i}\}$.
        \item 
        If $(0:_M I_k)=0$ for all $k=1,\dots,r$, and ${\bf I}^{\underline{s}} M\subseteq N$ for some $\underline{s}\in\N^r$, then for each $\fp\in\mathcal{A}_N^M({\bf I})$, the same result holds true for $v_\fp(M/{\bf I}^{\underline{n}} N)$, i.e., $v_\fp(M/{\bf I}^{\underline{n}} N)$ is asymptotically the minimum of finitely many linear functions in $\underline{n}$.
        \item 
        With the same hypotheses of $(2)$, given $\fp\in\mathcal{B}_{{\bf I}N}^M({\bf I})$ $($hence $\fp\in\mathcal{A}_N^M({\bf I})$$)$, the functions $v_\fp({\bf I}^{\underline{n}}M/{\bf I}^{\underline{n}+\underline{1}} N)$ and $v_\fp(M/{\bf I}^{\underline{n}+\underline{1}} N)$ coincide for all $\underline{n}\gg \underline{0}$.
    \end{enumerate}
\end{theorem}

The following result is a direct consequence of Theorem \ref{thm::asymbehaviour}.

\begin{corollary}\label{cor:asymp-v-num}
    With {\rm Setup~\ref{setup}}, the $v$-number $v({\bf I}^{\underline{n}}M/{\bf I}^{\underline{n}} N)$ eventually becomes either $\infty$, or the minimum of finitely many linear functions in $\underline{n}$. The same holds for the function $v(M/{\bf I}^{\underline{n}} N)$ under the additional conditions that $(0:_M I_k)=0$ for all $k=1,\ldots,r$, and ${\bf I}^{\underline{s}} M\subseteq N$ for some $\underline{s}\in\N^r$.
\end{corollary}

When $R = R_0[X_1,\ldots,X_d]$ is a (graded) polynomial ring over a Noetherian integral domain $R_0$, Corollary \ref{cor:asymp-v-num} yields that $v(R/{\bf I}^{\underline{n}})$ eventually is the minimum of finitely many linear functions in $\underline{n}$. Our next theorem shows that $v(R/{\bf I}^{\underline{n}})$ is, in fact, eventually a linear function in $\underline{n}$, where the leading coefficients are given by the initial degrees of $I_1,\ldots,I_r$. This result is surprising because $\reg(R/{\bf I}^{\underline{n}})$ is not always eventually linear even when $R$ is a polynomial ring over a field, as shown in \cite[Ex.~3.1]{BC17} by Bruns-Conca.

\begin{theorem}\label{thm:lin-v-function}
    Let $R = R_0[X_1,\ldots,X_d]$ be an $\N$-graded polynomial ring over a Noetherian integral domain $R_0$, and let $I_1,\ldots,I_r$ be nonzero homogeneous ideals such that $\indeg(I_i)\ge 1$ for at least one $i$. Then, the functions $v(R/{\bf I}^{\underline{n}})$, $v({\bf I}^{\underline{n}}/{\bf I}^{\underline{n}+\underline{1}})$ and $\indeg({\bf I}^{\underline{n}}/{\bf I}^{\underline{n}+\underline{1}})$ eventually become linear in $\underline{n}$ with the same leading coefficients given by $(d_1,\dots,d_r)$, where $d_i := \indeg(I_i)$ for $1\le i\le r$.
\end{theorem}

When $R_0$ is a field, in Theorem \ref{thm:lin-v-function}, the condition $\indeg(I_i)\ge 1$ is equivalent to that $I_i$ is a proper ideal of $R$.

Now we describe the contents of the article. In Section~\ref{sec::proofs}, we prove Theorems~\ref{thm::main} and \ref{thm:linearity}, which show the asymptotic behaviour of the initial degree and the (local) $v$-number of graded components in a finitely generated multigraded module. These results lead to the proofs of Theorems~\ref{thm::asymbehaviour} and \ref{thm:lin-v-function}. Finally, in Example~\ref{sec:Examples}, we provide some examples that complement Theorems~\ref{thm::asymbehaviour} and \ref{thm:lin-v-function}. In Examples~\ref{example-1} and \ref{example-2}, we see how the Vasconcelos invariant and the Castelnuovo–Mumford regularity can behave very differently, while Examples~\ref{example-2} and \ref{example-3} show some instances where the (local) $v$-number is not eventually linear, unlike the case for powers of a single ideal. Despite Theorem~\ref{thm:lin-v-function}, in Example~\ref{example-4}, we show that the linearity of local $v$-numbers cannot be expected even in the polynomial case. In Theorem~\ref{thm:lin-v-function}, the hypothesis on $R_0$ being an integral domain cannot be eluded, see Example~\ref{example-5}.

\section{Preliminaries on initial degrees and Gr\"obner bases}\label{sec:prelim}

Here, we provide some preliminaries on the initial degree and Gr\"obner bases that we use later in Section~\ref{sec::proofs}. We work under Setup~\ref{setup} and with the notations described in \ref{para:notations-r-tuple} and \ref{para:notations}. Note that if $M$ is nonzero, then $\indeg(M)$ is finite, since $M$ is finitely generated.

    \begin{lemma}\label{lemma:prop_indeg}
        Let $U,V$ and $W$ be $\Z$-graded $R$-modules. Then the following hold.
        \begin{enumerate}[\rm (1)]\label{lemma:indeg_firstprop}
            \item $\indeg(V(a))=\indeg(V)+a$, where $V(a)$ is the $\Z$-graded $R$-module whose $n$th graded component is $V_{n+a}$ for every $n\in\Z$.
            \item For a graded injective homomorphism $U\hookrightarrow V$, we have $\indeg(V)\le\indeg(U).$
            \item For a short exact sequence $0\to U\to V\to W\to 0$ of graded $R$-modules,
            $$\indeg(V) = \min \{\indeg(U),\indeg(W)\}.$$
            \item For a chain $0=V_0\subsetneq V_1\subsetneq\cdots\subsetneq V_q=V$ of graded submodules of $V$, $$\indeg(V)=\min\{\indeg(V_{i+1}/V_{i}):0\le i\le q-1\}.$$
        \end{enumerate}
    \end{lemma}
    \begin{proof}
        The first two statements follow from the definition of initial degree, while the third statement can be obtained from the observation that the $n$th graded component $V_n=0$ if and only if both $U_n=0$ and $W_n=0$. Finally, from the chain of graded submodules, one derives the following graded short exact sequences:
        $$0\to V_i\to V_{i+1}\to V_{i+1}/V_i\to 0\quad\text{for }0\le i\le q-1.$$
        The last statement is then a direct consequence of the third one.
    \end{proof}

We are interested in the behavior of the initial degree under Gr\"obner perturbation. The theory of Gr\"obner bases can also be developed for polynomial rings over general rings, with some precautions. Consider a polynomial ring $S := A[X_1,\dots,X_n]$ over a Noetherian ring $A$. As in the case of polynomial rings over a field, the ring $S$ is a free $A$-module whose canonical basis is given by the set of monomials of $S$, denoted by $\Mon(S)$. We consider a monomial order on $\Mon(S)$.
    
\begin{definition}\label{def:monomialorder}
    A \emph{monomial order}, or \emph{term order}, is a total order $<$ on $\Mon(S)$ satisfying the properties:
    \begin{enumerate}[(1)]
        \item $1< m$ for every $m\in \Mon(S)$, $m\neq 1$;
        \item if $m_1,m_2,m_3\in \Mon(S)$ with $m_1< m_2$, then $m_1m_3< m_2m_3$.
    \end{enumerate}
\end{definition}

The following is a well-known fact about monomial order, cf.~\cite[Thm.~1.4.6]{MR1287608}.

\begin{proposition}\label{prop:well-order}
    Every monomial order $<$ on $\Mon(S)$ is a well-ordering, that is, every non-empty subset of $\Mon(S)$ has a minimal element.
\end{proposition}

\begin{para}\label{para:initial-ideal}
    Fix a term order $<$ on $S$. Given an element $f\in S$, we can write this element as $f=a_1m_1+\dots+a_tm_t$, where $a_i\in A\smallsetminus\{0\}$ and $m_1,\dots,m_t\in \Mon(S)$ such that $m_t<\dots<m_2<m_1$. Then $m_1$ is called the \emph{initial monomial} of $f$; while the \emph{initial term} is $\ini_<(f):=a_1m_1$. For a given ideal $I$ of $S$, the \emph{initial ideal} of $I$ is a monomial ideal defined as $\ini_<(I):=(\ini_<(g):g\in I)$. By convention, $\ini_<(0):=0$. 
\end{para}

\begin{remark}
    In the case where $A$ is a field, it is equivalent to define the initial ideal of an ideal using either the initial monomials or the initial terms. However, in the general case when $A$ is not a field, these two notions give rise to two different ideals.
\end{remark}

\begin{remark}\label{thm:grobner}
    Let $I$ be an ideal of $S$. Then there exists a finite system of generators $\mathcal{G}=\{f_1,\dots,f_q\}$ of $I$ such that $\ini_<(I)$ is generated by $\{\ini_<(f_1),\dots,\ini_<(f_q)\}$, see, e.g., \cite[Cor.~4.1.17]{MR1287608}. The set $\mathcal{G}$ is called a \emph{Gr\"obner basis} for $I$.
\end{remark}
    
\begin{remark}\label{rmk:megapower}
    In view of Definition~\ref{def:monomialorder}.(2), Remark~\ref{thm:grobner} ensures that every term in $\ini_<(I)$ can be written as the initial term of some element in $I$.
\end{remark}

\begin{para}\label{para:termorder_freemodule}
    For our purposes, we consider the generalized theory of Gr\"obner bases, extending it from ideals to modules. Let $F$ be a finitely generated free $S$-module, and let $\mathcal{B}_F:=\{e_1,\dots,e_s\}$ be the standard basis of $F$. Then, as an $A$-module, $F$ is generated by the set $\{ Xe_i : X\in \Mon(S), 1\le i\le s\}$ of \emph{terms}. A natural way to construct a term order on $F$ is to fix a term order $<$ on $S$, and define
    $$X_ie_i<X_je_j \text{ if and only if } i<j, \text{ or }i=j \mbox{ and } X_i<X_j.$$
    A more general approach is to consider the symmetric algebra $R[e_1,\dots,e_s]$ and view $F$ as the $R$-submodule generated by $e_1,\dots,e_s$. Then any term order on the symmetric algebra restricts to a term order on the free module. Finally, as in the classical case, one can now consider initial term of every element of $F$. Given a subset $U$ of $F$, the initial submodule of $U$ is defined to be
    \begin{center}
        $\ini_<(U) :=$ the $A$-submodule of $F$ generated by $\{\ini_<(u):u\in U\}$.
    \end{center}
    It turns out that if $U$ is an $S$-submodule of $F$, then $\ini_<(U)$ is also an $S$-submodule of $F$. Indeed, $\ini_<(U)$ is always an $A$-submodule of $F$. Moreover, if $U$ is an $S$-submodule of $F$, then $\ini_<(U)$ is closed under multiplication by elements of $S$, since $\ini_<(U)$ is closed under multiplication by monomials and any polynomial in $S$ is an $A$-linear combination of monomials. Analogous results to Proposition~\ref{prop:well-order} and Theorem~\ref{thm:grobner} continue to hold in this setting. More details on this topic can be found in \cite[Sec.~1.4]{KR} and \cite[Chapter~3]{MR1287608}.
\end{para}
    
We are now ready to prove that the initial degree behaves nicely under Gr\"obner degeneration. The main idea of the following result is that by taking the initial module, no graded components can appear and/or disappear.
\begin{proposition}\label{prop:initialdegree_gbdegeneration}
    Let $\mathscr{F}$ be a $\Z^{r+1}$-graded finitely generated free $S$-module. Set a term order $<$ on $\mathscr{F}$. Let $\mathscr{U}$ be a nonzero multigraded submodule of $\mathscr{F}$. For every $\underline{n}\in\Z^{r}$, denote $\mathscr{U}_{(\underline{n},*)}:=\bigoplus_{l\in\Z} \mathscr{U}_{(\underline{n},l)}$ and $\mathscr{F}_{(\underline{n},*)}:=\bigoplus_{l\in\Z} \mathscr{F}_{(\underline{n},l)}$. Then
    $$\indeg \big(\mathscr{F}_{(\underline{n},*)}/\mathscr{U}_{(\underline{n},*)}\big) = \indeg\big(\mathscr{F}_{(\underline{n},*)}/\ini_<(\mathscr{U})_{(\underline{n},*)}\big) \;\text{ for every } \underline{n}\in\Z^{r}.$$
\end{proposition}

\begin{proof}
    We use the setup as discussed in \ref{para:termorder_freemodule}, i.e., $\rank(\mathscr{F})=s$, and $\mathscr{F}$ is generated as an $A$-module by the terms $\{ Xe_i : X\in \Mon(S), 1\le i\le s\}$. Fix $\underline{n}\in\Z^{r}$. In order to prove the desired equality, it is enough to show that, for every $k\in\Z$, $\mathscr{F}_{(\underline{n},k)}=\mathscr{U}_{(\underline{n},k)}$ if and only if $\mathscr{F}_{(\underline{n},k)}=\ini_<(\mathscr{U})_{(\underline{n},k)}$. We fix $k\in\Z$.

    Suppose $\mathscr{F}_{(\underline{n},k)}=\mathscr{U}_{(\underline{n},k)}$. Since $\mathscr{F}_{(\underline{n},k)}$ is generated, as an $A$-module, by the terms of $\mathscr{F}$ of degree $(\underline{n},k)$, it remains unaffected by taking the initial submodule. Therefore, one has the equality $\mathscr{F}_{(\underline{n},k)}=\ini_<(\mathscr{U}_{(\underline{n},k)})$. Now we notice that $\ini_<(\mathscr{U}_{(\underline{n},k)})\subset \ini_<(\mathscr{U})_{(\underline{n},k)}$. Indeed, given an element $u\in \mathscr{U}_{(\underline{n},k)}\subset \mathscr{U}$, its initial term $\ini_<(u)$ has degree $(\underline{n},k)$, hence $\ini_<(u)\in \ini_<(\mathscr{U})_{(\underline{n},k)}$.
    Thus $\mathscr{F}_{(\underline{n},k)} = \ini_<(\mathscr{U}_{(\underline{n},k)})\subseteq \ini_<(\mathscr{U})_{(\underline{n},k)}\subseteq \mathscr{F}_{(\underline{n},k)}$, which implies that $\mathscr{F}_{(\underline{n},k)} = \ini_<(\mathscr{U})_{(\underline{n},k)}$.
        
    For the reverse implication, we proceed by contradiction: Suppose that $\mathscr{F}_{(\underline{n},k)}=\ini_<(\mathscr{U})_{(\underline{n},k)}$, but $\mathscr{F}_{(\underline{n},k)}\neq \mathscr{U}_{(\underline{n},k)}$. Since $\mathscr{F}_{(\underline{n},k)}$ is generated as an $A$-module by the terms of $\mathscr{F}$ of degree $(\underline{n},k)$, there exists an element $Xe_i\in \mathscr{F}_{(\underline{n},k)} \smallsetminus\mathscr{U}_{(\underline{n},k)}$ for some $X\in \Mon(S)$ and $1 \le i \le s$. By the well-ordering of monomial orders, we can suppose that this is the smallest term in $\mathscr{F}_{(\underline{n},k)} \smallsetminus\mathscr{U}_{(\underline{n},k)}$.
        
    Since $\mathscr{F}_{(\underline{n},k)}=\ini_<(\mathscr{U})_{(\underline{n},k)}$, there exists a homogeneous element $u\in \mathscr{U}_{(\underline{n},k)}$ such that $\ini_<(u)=Xe_i$ by Remark~\ref{rmk:megapower}. The element $u$ can be written as $u=Xe_i+u_1$, where the monomials in $u_1$ are smaller than $Xe_i$. By the minimality of $Xe_i$, it follows that $u_1$ must lie in $\mathscr{U}_{(\underline{n},k)}$. Hence $Xe_i=u-u_1\in \mathscr{U}_{(\underline{n},k)}$, which is a contradiction.
\end{proof}

The local Vasconcelos invariant of a module $M$ at $\fp\in \Ass_R(M)$ can be interpreted as the initial degree of certain module depending on $M$ and $\fp$.

\begin{lemma}\cite[Lem.~1.5]{FG}\label{lem:$v$-num-indeg}
    With {\rm Setup~\ref{setup}}, let $\fp\in\Ass_R(M)$. Set $X_\fp := \{\fq \in \Ass_R(M) : \fp \subsetneq \fq \}$. Let $V=R$ if $X_\fp=\emptyset$, otherwise $V = \prod_{\fq\in X_\fp}\fq$. Then
    $$v_\fp(M) = \indeg\big(\ann_M(\fp)/\ann_M(\fp)\cap\Gamma_V(M)\big).$$
\end{lemma}

\section{Proof of the results}\label{sec::proofs}

The following result is a generalization of \cite[Thm.~2.8]{FG}. In the proof of asymptotic behaviour of $\indeg\left(  \mathscr{L}_{(\underline{n},*)} \right)$, we use an argument similar to \cite[Proof of Thm.~8.3.7]{BCRV}.

\begin{theorem}\label{thm::main}
    Let $ T = R_0 [x_1,\dots,x_d,y_{1,1},\ldots,y_{1,a_1},\ldots,y_{r,1},\dots,y_{r,a_r}] $ be a $\mathbb{Z}^{r+1}$-graded ring over a commutative Noetherian ring $R_0$, where $ \deg(x_{i}) = (\underline{0},f_{i})$ for $1\le i\le d$ and $\deg(y_{i,j}) = (\underline{e}_i,d_{i,j})$ for $ 1 \le i \le r$, $1 \le j \le a_i$. Assume that $f_i\ge 0$ for $1\le i\le d$. Let $\mathscr{L}$ be a finitely generated $\mathbb{Z}^{r+1}$-graded $T$-module. Set $ R := R_0[x_1,\dots,x_d] $, where $\deg(x_i) = f_i$ for $1\le i\le d$. Denote $\mathscr{L}_{(\underline{n},*)} := \bigoplus_{l\in\mathbb{Z}}  \mathscr{L}_{(\underline{n},l)}$ for each $\underline{n} \in \mathbb{Z}^r$.
	
    Note that $R$ is an $\N$-graded ring, and $ \mathscr{L}_{(\underline{n},*)} $ is a finitely generated $\mathbb{Z}$-graded $R$-module for each $\underline{n} \in \mathbb{Z}^r$. Moreover, the set $\Ass_R(\mathscr{L}_{(\underline{n},*)})$ stabilizes to a set, say $\mathcal{A}_{\mathscr{L}}$, for all $\underline{n}\gg \underline{0}$. It follows that $\mathscr{L}_{(\underline{n},*)} = 0$ for all $\underline{n}\gg \underline{0}$, or $\mathscr{L}_{(\underline{n},*)} \neq 0$ for all $\underline{n}\gg \underline{0}$. Assume the second case. Suppose $F(\underline{n})=\indeg\left(  \mathscr{L}_{(\underline{n},*)} \right)$, or $F(\underline{n})=v_{\fp}\left(  \mathscr{L}_{(\underline{n},*)} \right)$ for $\fp \in \mathcal{A}_{\mathscr{L}}$, or $F(\underline{n})=v\left(\mathscr{L}_{(\underline{n},*)} \right)$ for all $\underline{n}\in\bz^r$.
        
    Then, there exist $\underline{\omega}_1,\ldots,\underline{\omega}_s\in\Z^r$ and $c_1,\ldots,c_s\in\mathbb{Z}$, depending on $F$, such that
    $$F(\underline{n})=\min\{\underline{\omega}_j\cdot\underline{n}+c_j : 1\le j\le s\} \ \mbox{ for all } \underline{n}\gg \underline{0},$$
    where the $i$th component $\omega_{ji}$ of the coefficient vector $\underline{\omega}_j$ lies in $\{d_{i,1},\dots,d_{i,a_i}\}$ for $1\le i\le r$ and $1\le j\le s$. Recall that $\underline{\omega}_j \cdot\underline{n} = \omega_{j1} n_1 + \cdots + \omega_{jr} n_r$ for $\underline{n}\in\bz^r$.
\end{theorem}

\begin{proof}
    By writing $T=R[y_{1,1},\ldots,y_{1,a_1},\ldots,y_{r,1},\dots,y_{r,a_r}]$ with $\deg(y_{i,j}) = \underline{e_i}$ for $ 1 \le i \le r$, $1 \le j \le a_i$, we can realize $T$ as a Noetherian standard $\N^r$-graded ring over $T_{\underline{0}}=R$. Thus $\mathscr{L}=\bigoplus_{\underline{n}\in\Z^r}\mathscr{L}_{(\underline{n},*)}$ becomes a finitely generated $\Z^r$-graded $T$-module. So, by \cite[Thm.~3.4.(i)]{West04}, the set $\Ass_R(\mathscr{L}_{(\underline{n},*)})$ stabilizes to a set, say $\mathcal{A}_{\mathscr{L}}$, for $\underline{n}\gg \underline{0}$. If $\mathcal{A}_{\mathscr{L}}$ is an empty set, then $\mathscr{L}_{(\underline{n},*)}=0$ for all $\underline{n}\gg \underline{0}$. In the second case, assume that $\mathcal{A}_{\mathscr{L}}\not=\emptyset$. In this case, $\mathscr{L}_{(\underline{n},*)}\not=0$ for all $\underline{n}\gg \underline{0}$.

    (1) We first prove that $\indeg\left(  \mathscr{L}_{(\underline{n},*)} \right)$ is asymptotically the minimum of finitely many linear functions. Consider the polynomial ring
    $$ \mathcal{T} := R[Y_{1,1},\ldots,Y_{1,a_1},\ldots,Y_{r,1},\dots,Y_{r,a_r}],$$
    where $\deg(f) = (\underline{0},\deg_R(f))$ for $f\in R$ and $\deg(Y_{i,j}) = (\underline{e}_i,d_{i,j})$. Then $\mathscr{L}$ can be regarded as a $\mathcal{T}$-module via the natural ring homomorphism $\mathcal{T} \to T$. We start by presenting $\mathscr{L}$ as a quotient $\mathscr{F}/\mathscr{U}$, where $\mathscr{F}$ is a $\Z^{r+1}$-graded free $\mathcal{T}$-module, and $\mathscr{U}$ is a multigraded submodule of $\mathscr{F}$. Then, by taking any term order $<$ on $\mathscr{F}$, as explained in \ref{para:termorder_freemodule}, consider the initial submodule $\init(\mathscr{U})$. From Proposition~\ref{prop:initialdegree_gbdegeneration}, it follows that $\indeg\left(\mathscr{L}_{(\underline{n},*)}\right)=\indeg\left(\mathscr{F}_{(\underline{n},*)}/(\init (\mathscr{U}))_{(\underline{n},*)}\right)$. Next consider a chain of multigraded submodules
    $$0=\mathscr{M}^0\subsetneq \mathscr{M}^1\subsetneq\dots\subsetneq\mathscr{M}^s=\mathscr{F}/(\init(\mathscr{U}))$$
    in such a way that any consecutive quotient $\mathscr{M}^j/\mathscr{M}^{j-1}$ is isomorphic to a quotient of $\mathcal{T}$ by a monomial prime ideal (up to a degree shift). In particular, by restricting the chain to the $(\underline{n},*)$ graded component and by applying Lemma~\ref{lemma:prop_indeg}.(4), we obtain
    \[
    \indeg\left(\mathscr{L}_{(\underline{n},*)}\right)=\min\left\{\indeg \left(\mathscr{M}^{j}_{(\underline{n},*)}/\mathscr{M}^{j-1}_{(\underline{n},*)}\right): 1\le j\le s \right\}.
    \]
    Hence it is enough to show that $\indeg \big(\mathscr{M}^{j}_{(\underline{n},*)}/\mathscr{M}^{j-1}_{(\underline{n},*)}\big)$ is asymptotically a linear function in $\underline{n}$ whose leading coefficients are taken from the degrees $d_{i,j}$. Thus, we reduce to the case that $\mathscr{L}=(\mathcal{T}/\mathcal{J})(-\underline{u},-b)$ for some $\underline{u}\in\Z^r$ and $b\in \Z$, where $\mathcal{J}=J_0\mathcal{T}+(Y_{i,j}:Y_{i,j}\notin V)$ for some prime ideal $J_0$ of $R$ and for some subset $V$ of the set of the variables $\{Y_{1,1},\ldots,Y_{1,a_1},\ldots,Y_{r,1},\dots,Y_{r,a_r}\}$. Since $\mathscr{L}_{(\underline{n},*)}\not=0$ for all $\underline{n}\gg \underline{0}$, the intersection $V\cap\{Y_{i,1},\dots Y_{i,a_i}\}$ is not an empty set for every $1\le i \le r$. Set $w_i:=\min\{d_{i,j}: 1 \le j\le a_i, \,Y_{i,j}\in V\}$ for $i=1,\dots,r$, and $\underline{w}:=(w_1,\dots,w_r)$. Hence, since $\mathscr{L}=(\mathcal{T}/\mathcal{J})(-\underline{u},-b)$, it follows that
    \[ \indeg\left(\mathscr{L}_{(\underline{n},*)}\right)=\underline{w}\cdot(\underline{n}-\underline{u})+b=\underline{w}\cdot\underline{n}+\Tilde{b},
    \]
    where $\Tilde{b} := b-\underline{w} \cdot \underline{u}$. Note that here we need the condition that $f_i\ge0$ for $1\le i\le d$.

    (2) We now prove that $v_{\fp}\left(  \mathscr{L}_{(\underline{n},*)} \right)$ for $\fp \in \mathcal{A}_{\mathscr{L}}$ is asymptotically the minimum of finitely many linear functions. Given $\fp\in\mathcal{A}_{\mathscr{L}}$, denote $X_\fp := \{\fq \in \mathcal{A}_{\mathscr{L}} : \fp \subsetneq \fq \}$. Set $V := \prod_{\fq\in X_\fp}\fq$ (in the critical case $X_\fp=\emptyset$, set $V:=R$). Let $\mathscr{M} = \ann_{\mathscr{L}}(\fp)/\ann_{\mathscr{L}}(\fp)\cap\Gamma_V(\mathscr{L})$. Since $ T $ is Noetherian, and $\mathscr{L}$ is finitely generated, the quotient $\mathscr{M}$ is also a finitely generated $\mathbb{Z}^{r+1}$-graded $T$-module, whose grading is induced by that of $\mathscr{L}$. In particular, $ \mathscr{M}_{(\underline{n},*)} = \bigoplus_{l\in\mathbb{Z}}  \mathscr{M}_{(\underline{n},l)}$ is same as $\ann_{\mathscr{L}_{(\underline{n},*)}}(\fp)/\ann_{\mathscr{L}_{(\underline{n},*)}}(\fp)\cap\Gamma_V({\mathscr{L}_{(\underline{n},*)}})$ for every $\underline{n}\in\Z^{r}$. Therefore, Lemma~\ref{lem:$v$-num-indeg} together with (1) yields the desired result for $v_{\fp}\left( \mathscr{L}_{(\underline{n},*)} \right)$.
	
	(3) Note that $\mathcal{A}_{\mathscr{L}}$ is a non-empty finite set. By the definition of Vasconcelos invariant, $v\left( \mathscr{L}_{(\underline{n},*)} \right) = \inf\{ v_{\fp}\left( \mathscr{L}_{(\underline{n},*)} \right) : \fp\in\mathcal{A}_{\mathscr{L}} \}$ for all $\underline{n}\gg \underline{0}$. Hence, by (2), $v\left(  \mathscr{L}_{(n,*)} \right) $ is eventually the minimum of finitely many linear functions whose leading coefficients are taken from the degrees $d_{i,j}$.
\end{proof}

\begin{remark}
    (1) The condition $f_i\ge0$ is a strong condition that makes sure that Theorem~\ref{thm::main} holds true. Indeed, suppose that $f_i<0$ for some $i=1,\dots,d$, and $x_i$ is $\mathscr{L}$-regular. Then, by taking $0\neq \ell\in\mathscr{L}_{(\underline{n},*)}$, the element $x_i^k\cdot\ell$ is nonzero in $\mathscr{L}_{(\underline{n},*)}$ for every $k\in\N$, which implies that $\indeg(\mathscr{L}_{(\underline{n},*)}) = -\infty$.

    (2) In Theorem~\ref{thm::main}, either $\mathscr{L}_{(\underline{n},*)} = 0$ for all $\underline{n}\gg \underline{0}$, or $\mathscr{L}_{(\underline{n},*)} \neq 0$ for all $\underline{n}\gg \underline{0}$. In the second case, $\indeg(\mathscr{L}_{(\underline{n},*)})$ is finite for all $\underline{n}\gg \underline{0}$, since each $\mathscr{L}_{(\underline{n},*)} $ is a finitely generated $\mathbb{Z}$-graded module over $R = R_0[x_1,\dots,x_d]$, where $\deg(x_i) = f_i\ge 0$.
\end{remark}

Under some additional conditions, the functions $\indeg\left(\mathscr{L}_{(\underline{n},*)}\right)$ and $v\left(\mathscr{L}_{(\underline{n},*)} \right)$ in Theorem~\ref{thm::main} are eventually linear in $\underline{n}$, as shown below.

\begin{theorem}\label{thm:linearity}
    With the hypotheses as in {\rm Theorem~\ref{thm::main}}, without loss of generality, assume that $d_{i,1}\le d_{i,2}\le \cdots\le d_{i,a_i}$ for $1\le i\le r$. Let $y_{1,1}\cdots y_{r,1}\notin\sqrt{\ann \mathscr{L}}$. Then, the functions $\indeg\left(\mathscr{L}_{(\underline{n},*)}\right)$ and $v\left(\mathscr{L}_{(\underline{n},*)} \right)$ become linear for all $\underline{n}\gg \underline{0}$ with the same leading coefficients given by $\underline{\delta}:=(d_{1,1},d_{2,1},\dots,d_{r,1})$.
\end{theorem}

\begin{proof}
    Set ${\bf y}^{\underline{n}}:=y_{1,1}^{n_1}\cdots y_{r,1}^{n_r}$ for $\underline{n}\in\N^r$. Then $\deg({\bf y}^{\underline{n}}) = (\underline{n},\,\underline{\delta}\cdot\underline{n})$ for all $\underline{n}\in\N^r$. Suppose $ \mathscr{L} = \bigoplus_{\underline{n}\in\Z^r} \mathscr{L}_{(\underline{n},*)} $ is generated by homogeneous elements of degree $\le\underline{m}$. We first prove the following claims:
    \begin{enumerate}[]
        \item {\it Claim 1}. There exists $\ell\in\N$ such that $(0:_{\mathscr{L}}{\bf y}^{\underline{n}})=(0:_{\mathscr{L}}{\bf y}^{
        \ell\cdot\underline{1}})$ for every $\underline{n}\ge \ell\cdot\underline{1}$.
        \item {\it Claim 2}. For every $n\in\N$, $\big(0:_{\mathscr{L}_{(\underline{m},*)}}{\bf y}^{n\cdot\underline{1}}\big)$ is a proper submodule of $\mathscr{L}_{(\underline{m},*)}$.
        \item {\it Claim 3}.
        For every $\underline{n}\in\Z^r$ and $\underline{\nu}\in\N^r$, one has
        \begin{equation*}
            \indeg\left(\mathscr{L}_{(\underline{n},*)}\right) \le \indeg\left( \dfrac{\mathscr{L}_{(\underline{n}-\underline{\nu},*)}}{\big(0:_{\mathscr{L}_{(\underline{n}-\underline{\nu},*)}} {\bf y}^{\underline{\nu}}\big)}  \right) + \underline{\delta}\cdot\underline{\nu}.
        \end{equation*}
        The same inequality holds for the $v$-numbers and the local $v$-numbers at every associated prime of the quotient $R$-module in the right hand side.
    \end{enumerate}

    {\it Proof of Claim 1.} Since the module $\mathscr{L}$ is Noetherian, the chain of submodules
    $$(0:_{\mathscr{L}} {\bf y}^{\underline{1}}) \subseteq (0:_{\mathscr{L}} {\bf y}^{2\cdot\underline{1}}) \subseteq (0:_{\mathscr{L}} {\bf y}^{3\cdot\underline{1}}) \subseteq \cdots $$
    stabilizes. So there exists $\ell\ge 1$ such that $(0:_{\mathscr{L}} {\bf y}^{n\cdot\underline{1}}) = (0:_{\mathscr{L}} {\bf y}^{\ell\cdot\underline{1}})$ for every $n\ge \ell$. Fix $\underline{n}\in\N^r$ such that $\underline{n}\ge \ell\cdot\underline{1}$. Set $\alpha:=\max\{n_i:1\le i\le r\}$. Then one has $\ell\cdot\underline{1}\le\underline{n}\le \alpha\cdot\underline{1}$, which implies that
    \[
    (0:_{\mathscr{L}} {\bf y}^{\ell\cdot\underline{1}}) \subseteq (0:_{\mathscr{L}} {\bf y}^{\underline{n}}) \subseteq (0:_{\mathscr{L}} {\bf y}^{\alpha\cdot\underline{1}}).
    \]
    Since the submodules on both sides coincide by the construction of $\ell$, they must all coincide to $(0:_{\mathscr{L}} {\bf y}^{\ell\cdot\underline{1}})$. This proves Claim 1.
    
    {\it Proof of Claim 2.} If possible, let $(0:_{\mathscr{L}_{(\underline{m},*)}}{\bf y}^{n\cdot\underline{1}})=\mathscr{L}_{(\underline{m},*)}$. Then ${\bf y}^{n\cdot\underline{1}}\mathscr{L}_{(\underline{m},*)}=0$. Since $\mathscr{L}$ is finitely generated in degrees $\le \underline{m}$, it follows that ${\bf y}^{\underline{1}} = y_{1,1}\cdots y_{r,1}\in\sqrt{\ann \mathscr{L}}$, which is a contradiction. So $(0:_{\mathscr{L}_{(\underline{m},*)}}{\bf y}^{n\cdot\underline{1}}) \subsetneqq\mathscr{L}_{(\underline{m},*)}$.

    {\it Proof of Claim 3.} Fix $\underline{n}\in\Z^r$ and $\underline{\nu}\in\N^r$. Consider the $T$-module homomorphism $\mathscr{L}\to\mathscr{L}$ given by multiplication with ${\bf y}^{\underline{\nu}}$. Since $\deg({\bf y}^{\underline{\nu}}) = (\underline{\nu},\,\underline{\delta}\cdot\underline{\nu})$, it induces an injective graded $R$-module homomorphism
    \begin{equation}\label{inj-map-colon-L}
    \dfrac{\mathscr{L}_{(\underline{n}-\underline{\nu},*)}}{\big(0:_{\mathscr{L}_{(\underline{n}-\underline{\nu},*)}} {\bf y}^{\underline{\nu}}\big)}(-(\underline{\delta}\cdot\underline{\nu}))
    \stackrel{{\bf y}^{\underline{\nu}}}{\longhookrightarrow} \mathscr{L}_{(\underline{n},*)}.    
    \end{equation}    
    Here $M(-m)$ denotes the graded $R$-module with $M_{n-m}$ as its $n$th graded component. By the definition of $v$-numbers, $v_\fp(M(-m)) = v_\fp(M)+m$ for all $\fp\in\Ass_R(M)$. Claim 3 now follows from \eqref{inj-map-colon-L} using Lemma~\ref{lemma:prop_indeg} (1)-(2) and \cite[Prop.~2.5]{FG}.

    Set $\underline{n}_0:=\underline{m}+\ell\cdot\underline{1}$. Combining the three claims above, for every $\underline{n}\ge\underline{n}_0$, considering $\underline{\nu}=\underline{n}-\underline{m}$ in Claim 3, one obtains that
    \begin{align}
        \indeg\left(\mathscr{L}_{(\underline{n},*)}\right) &\le \indeg\left( \dfrac{\mathscr{L}_{(\underline{m},*)}}{\big(0:_{\mathscr{L}_{(\underline{m},*)}} {\bf y}^{\underline{n}-\underline{m}}\big)}  \right) + (\underline{\delta}\cdot(\underline{n}-\underline{m})) \label{ineq-indeg-quotien} \\
        &=\underline{\delta}\cdot\underline{n}+\indeg\left( \dfrac{\mathscr{L}_{(\underline{m},*)}}{\big(0:_{\mathscr{L}_{(\underline{m},*)}} {\bf y}^{\ell\cdot\underline{1}}\big)}  \right) - (\underline{\delta}\cdot\underline{m})<\infty.\label{eq-indeg-quotient}
    \end{align}    
    Thus, there exists $c\in\Z$ such that
    \begin{equation}\label{indeg-L-bound}
        \indeg\left(\mathscr{L}_{(\underline{n},*)}\right)\le\underline{\delta}\cdot\underline{n}+c \mbox{ for all } \underline{n} \ge \underline{n}_0.
    \end{equation}
    On the other hand, in Theorem~\ref{thm::main}, it is shown that there exist $\underline{\omega}_1,\ldots,\underline{\omega}_s\in\Z^r$ and $c_1,\ldots,c_s\in\mathbb{Z}$ such that
    \begin{equation}\label{indeg-L-equality}
        \indeg\left(\mathscr{L}_{(\underline{n},*)}\right)=\min\{\underline{\omega}_j\cdot\underline{n}+c_j : 1\le j\le s\} \ \mbox{ for all } \underline{n}\gg \underline{0},
    \end{equation}
    where the $i$th component $\omega_{ji}$ of the coefficient vector $\underline{\omega}_j$ lies in $\{d_{i,1},\dots,d_{i,a_i}\}$ for $1\le i\le r$ and $1\le j\le s$. In particular, by the given hypothesis, $\underline{\omega}_j\ge \underline{\delta}$ for $1\le j\le s$, which yields that $\underline{\omega}_j\cdot\underline{n}\ge \underline{\delta}\cdot\underline{n}$ for all $\underline{n}\in\N^r$. Thus, combining \eqref{indeg-L-bound} and \eqref{indeg-L-equality}, there exists $b\in\Z$ such that
    \begin{equation}\label{eqn::squeeze}
        \underline{\delta}\cdot\underline{n}+b\le \indeg\left(\mathscr{L}_{(\underline{n},*)}\right)\le \underline{\delta}\cdot\underline{n}+c
        \quad\text{for all }\underline{n}\gg \underline{0}.
    \end{equation}
    Hence, for every fixed $\underline{\nu}\gg \underline{0}$, one has that
    \begin{equation}\label{eqn:squeeze-m}
        m(\underline{\delta}\cdot\underline{\nu})+b\le \indeg\left(\mathscr{L}_{(m\underline{\nu},*)}\right)\le m(\underline{\delta}\cdot\underline{\nu})+c
            \quad\text{for all } m \gg 0.
    \end{equation}
    On the other hand, for every fixed $\underline{\nu}\gg \underline{0}$, by \eqref{indeg-L-equality}, the function $\indeg\left(\mathscr{L}_{(m\underline{\nu},*)}\right)$ is linear in $m$ for all $m\gg0$, in fact, there exists some $j\in\{1,\dots,s\}$ such that $\indeg\left(\mathscr{L}_{(m\underline{\nu},*)}\right)=m(\underline{w}_j\cdot\underline{\nu})+c_j$ for all $m\gg0$. In view of \eqref{eqn:squeeze-m}, the leading coefficient must be the same as $\underline{\delta}\cdot\underline{\nu}$. So $\underline{w}_j\cdot\underline{\nu} = \underline{\delta}\cdot\underline{\nu}$ for all $\underline{\nu}\gg \underline{0}$. Since $\underline{w}_j\ge\underline{\delta}$, it follows that $\underline{w}_j=\underline{\delta}$. Thus there exists $j\in\{1,\dots,s\}$ such that $\underline{w}_j=\underline{\delta}$. Set $a:=\min\{c_l : 1\le l \le s, \underline{w}_l = \underline{\delta}\}$. Then, by \eqref{indeg-L-equality},
    \[
    \indeg\left(\mathscr{L}_{(\underline{n},*)}\right)=\underline{\delta}\cdot\underline{n}+a \mbox{ for all } \underline{n} \gg \underline{0}.
    \]   

    Similar inequalities as in \eqref{ineq-indeg-quotien} and \eqref{eq-indeg-quotient} for $v$-numbers yield that
    \[
        v(\mathscr{L}_{(\underline{n},*)})\le \underline{\delta}\cdot\underline{n}+e \mbox{ for all } \underline{n} \gg \underline{0},
    \]
    where $e\in\Z$. These are the inequalities like \eqref{indeg-L-bound}. Now, arguing in the same manner as for the function $\indeg(\mathscr{L}_{(\underline{n},*)})$, one obtains that $v(\mathscr{L}_{(\underline{n},*)})$ is eventually linear in $\underline{n}$ with the leading coefficients given by $\underline{\delta}$.
\end{proof}

\begin{remark}\label{rmk:local-v-lin}
    In the proof of Theorem~\ref{thm:linearity}, denote the quotient $R$-module considered in \eqref{eq-indeg-quotient} by $V$, i.e., $V:=\mathscr{L}_{(\underline{m},*)}/\big(0:_{\mathscr{L}_{(\underline{m},*)}} {\bf y}^{\ell\cdot\underline{1}}\big)$. Then, the injective homomorphisms in \eqref{inj-map-colon-L} yield that $\Ass_R(V)\subseteq \Ass_R(\mathscr{L}_{(\underline{n},*)})$ for all $\underline{n}\ge\underline{n}_0$. Hence, for every fixed $\fp\in\Ass_R(V)$,  following the same steps as \eqref{ineq-indeg-quotien} and \eqref{eq-indeg-quotient},
    \[
        v_\fp(\mathscr{L}_{(\underline{n},*)})\le \underline{\delta}\cdot\underline{n}+h \mbox{ for all } \underline{n} \gg \underline{0},
    \]
    where $h\in\Z$. These inequalities are obtained under the same considerations as \eqref{indeg-L-bound}. Now, arguing in the same manner, for $\fp\in\Ass_R(V)$, one sees that $v_\fp(\mathscr{L}_{(\underline{n},*)})$ is eventually linear in $\underline{n}$ with the leading coefficients given by $\underline{\delta}$.
\end{remark}

We are now in a position to prove the main theorems.

\begin{proof}[Proof of Theorem~\ref{thm::asymbehaviour}]
    Suppose $ R = R_0[x_1,\dots,x_d] $, where $\deg(x_i) = f_i$ for $1\le i\le d$. Let $I_i$ be generated by homogeneous elements $y_{i,1},\ldots,y_{i,a_i}$, where $\deg(y_{i,j}) = d_{i,j}$ for $1 \le j \le a_i$. We consider the Rees ring $\mathscr{R}=\mathscr{R}(I_1,\dots,I_r)$ with $\N^{r+1}$-graded structure given by $ \mathscr{R}_{(\underline{n},m)} = ({\bf I}^{\underline{n}})_m $ for all $(\underline{n},m)\in \N^{r+1}$. Thus, $\mathscr{R}$ can be identified with the multigraded ring $T$ as described in Theorem~\ref{thm::main}.

    (1) Let $\mathscr{R}(I_1,\dots,I_r;M)$ denote the Rees module of $M$ with respect to the ideals $I_1,\dots,I_r$. Set $\mathscr{L} := \mathscr{R}(I_1,\dots,I_r;M)/\mathscr{R}(I_1,\dots,I_r;N)$,
    where the grading is given by $\mathscr{L}_{(\underline{n},l)}:=({\bf I}^{\underline{n}}M/{\bf I}^{\underline{n}} N)_l$. Clearly, $\mathscr{L}$ is a finitely generated $\mathbb{Z}^{r+1}$-graded $\mathscr{R}$-module. Hence Theorem~\ref{thm::asymbehaviour}.(1) is a direct consequence of Theorem~\ref{thm::main}.
        
    (2) Let $\fp\in \mathcal{A}_N^M({\bf I})$. Set $X_\fp := \{\fq \in \mathcal{A}_N^M({\bf I}) : \fp \subsetneq \fq \}$. Let $V=R$ if $X_\fp=\emptyset$, otherwise $V = \prod_{\fq\in X_\fp}\fq$. Consider $\mathscr{G} := \mathscr{R}(I_1,\dots,I_r;M)/\mathscr{R}(I_1,\dots,I_r;{\bf I}N)$, which is a finitely generated $\Z^{r+1}$-graded $\mathscr{R}$-module. We now consider $\mathscr{L} := \ann_{\mathscr{G}}(\fp)/\ann_{\mathscr{G}}(\fp)\cap\Gamma_V(\mathscr{G})$. This is also a finitely generated $\Z^{r+1}$-graded $\mathscr{R}$-module, where the grading is induced by the one in $\mathscr{G}$. Using the notations as in Theorem~\ref{thm::main}, observe that
    \[
        \mathscr{L}_{(\underline{n},*)}=\dfrac{\ann_{{\bf I}^{\underline{n}}M/{\bf I}^{\underline{n}+\underline{1}}N}(\mathfrak{p})}{\ann_{{\bf I}^{\underline{n}}M/{\bf I}^{\underline{n}+\underline{1}}N}(\mathfrak{p})\cap\Gamma_{V}({\bf I}^{\underline{n}}M/{\bf I}^{\underline{n}+\underline{1}}N)} \mbox{ for all } \underline{n}\in\N^r.
    \]    
    Since ${\bf I}^{\underline{s}} M\subseteq N$, it follows that ${\bf I}\subseteq \fp$. Therefore
    \[
    ({\bf I}^{\underline{n}+\underline{1}}N:_M \fp)\subseteq ({\bf I}^{\underline{n}+\underline{1}}M:_M \fp)\subseteq ({\bf I}^{\underline{n}+\underline{1}}M:_M {\bf I}) = {\bf I}^{\underline{n}} M \mbox{ for all }\underline{n}\gg\underline{0},
    \]
    where the last equality is obtained by \cite[Lem.~1.3.(ii)]{KSharp}. Hence, a similar proof as that of \cite[Lem.~2.13]{FG} yields that
    \[
        \mathscr{L}_{(\underline{n},*)}=\dfrac{\ann_{M/{\bf I}^{\underline{n}+\underline{1}}N}(\mathfrak{p})}{\ann_{M/{\bf I}^{\underline{n}+\underline{1}}N}(\mathfrak{p})\cap\Gamma_{V}(M/{\bf I}^{\underline{n}+\underline{1}}N)}  \mbox{ for all } \underline{n}\gg \underline{0}.
    \]
    By Lemma~\ref{lem:$v$-num-indeg}, one has the equality $v_\fp(M/{\bf I}^{\underline{n}+\underline{1}}N)=\indeg (\mathscr{L}_{(\underline{n},*)})$ for all $\underline{n}\gg\underline{0}$. Theorem~\ref{thm::asymbehaviour}.(2) is now a consequence of Theorem~\ref{thm::main}. 
        
    (3) Given $\fp\in\mathcal{B}_{{\bf I}N}^M({\bf I})$. Then $\fp\in \mathcal{A}_N^M({\bf I})$. Following the notations as in the proof of (2), the functions $v_\fp({\bf I}^{\underline{n}}M/{\bf I}^{\underline{n}+\underline{1}} N)$ and $v_\fp(M/{\bf I}^{\underline{n}+\underline{1}} N)$ coincide for all $\underline{n}\gg \underline{0}$ since they both are asymptotically equal to $\indeg (\mathscr{L}_{(\underline{n},*)})$ by Lemma~\ref{lem:$v$-num-indeg}.
\end{proof}

\begin{remark}
    Using the assumptions and notations of Theorem~\ref{thm::asymbehaviour}, it is clear that given $\fp\in\mathcal{B}_{{\bf I}N}^M({\bf I})\subseteq \mathcal{A}_N^M({\bf I})$, the functions $v_\fp({\bf I}^{\underline{n}}M/{\bf I}^{\underline{n}+\underline{1}} N)$ and $v_\fp(M/{\bf I}^{\underline{n}+\underline{1}} N)$ coincide as long as $({\bf I}^{\underline{n}+\underline{1}}N:_M \fp)\subseteq {\bf I}^{\underline{n}} M$.
\end{remark}

\begin{proof}[Proof of Theorem~\ref{thm:lin-v-function}]
    Suppose $ R = R_0[X_1,\dots,X_d] $, where $\deg(X_i) = f_i$ for $1\le i\le d$. Here $f_i\ge 0$ for $1\le i\le d$. Let $I_i$ be generated by homogeneous elements $y_{i,1},\ldots,y_{i,a_i}$, where $\deg(y_{i,j}) = d_{i,j}$ for $1 \le j \le a_i$. Without loss of generality, we may assume that $d_{i,1}\le d_{i,2}\le \cdots\le d_{i,a_i}$ for $1\le i\le r$. Then $\indeg(I_i) = d_{i,1}$ for $1\le i\le r$. Consider the Rees ring $\mathscr{R}=\mathscr{R}(I_1,\dots,I_r)$, which can be identified with the multigraded ring $T$ as described in Theorem~\ref{thm::main}. Set $\mathscr{L}:=\mathscr{R}(I_1,\dots,I_r)/{\bf I}\mathscr{R}(I_1,\dots,I_r)$. Then $\mathscr{L}$ is a finitely generated $\N^{r+1}$-graded $\mathscr{R}$-module. Now, we follow the notations as in Theorems~\ref{thm::main} and \ref{thm:linearity}.
    
    We prove that the initial degree and the global $v$-number of $\mathscr{L}_{(\underline{n},*)}={\bf I}^{\underline{n}}/{\bf I}^{\underline{n}+\underline{1}}$ are eventually linear in $\underline{n}$ with the same leading coefficients given by $\underline{\delta}$. For this, in view of Theorem~\ref{thm:linearity}, it is enough to show that ${\bf y}:=y_{1,1}\cdots y_{r,1}\notin\sqrt{\ann_{\mathscr{R}} (\mathscr{L})}$. If possible, let ${\bf y} \in \sqrt{\ann_{\mathscr{R}} (\mathscr{L})}$. Then ${\bf y}^s\mathscr{L}=0$ for some $s\ge 1$. Since $\mathscr{L}_{(\underline{n},*)} = {\bf I}^{\underline{n}}/{\bf I}^{\underline{n}+\underline{1}}$ and $\deg({\bf y}) = (\underline{1}, \,\underline{\delta}\cdot\underline{1})$,
    it follows that ${\bf y}^s{\bf I}^{\underline{n}} \subseteq {\bf I}^{\underline{n}+(s+1)\cdot\underline{1}}$ for all $\underline{n}\in\N^r$. Denote $|\underline{\delta}| := \underline{\delta}\cdot\underline{1}$. As $R$ is an integral domain, ${\bf I}^{\underline{n}}\neq 0$, in addition $\indeg({\bf y}^s{\bf I}^{\underline{n}}) = s |\underline{\delta}|+\indeg({\bf I}^{\underline{n}})$ and $\indeg({\bf I}^{\underline{n}+(s+1)\cdot\underline{1}}) = (s+1)|\underline{\delta}|+\indeg({\bf I}^{\underline{n}})$.
    Thus
    \begin{align*}
        s|\underline{\delta}|+\indeg({\bf I}^{\underline{n}}) &= \indeg({\bf y}^s{\bf I}^{\underline{n}}) \\
        &\ge \indeg({\bf I}^{\underline{n}+(s+1)\cdot\underline{1}}) \quad \mbox{[by Lemma~\ref{lemma:indeg_firstprop} as ${\bf y}^s{\bf I}^{\underline{n}} \subseteq {\bf I}^{\underline{n}+(s+1)\cdot\underline{1}}$]} \\
        &= (s+1)|\underline{\delta}|+\indeg({\bf I}^{\underline{n}}),
    \end{align*}
    which is a contradiction as $|\underline{\delta}|\ge 1$. So ${\bf y} \notin \sqrt{\ann\mathscr{L}}$. This proves the result for the functions $v({\bf I}^{\underline{n}}/{\bf I}^{\underline{n}+\underline{1}})$ and $\indeg({\bf I}^{\underline{n}}/{\bf I}^{\underline{n}+\underline{1}})$.

    Note that $(0:_R I_i)=0$ for $1\le i\le r$. So, by Theorem~\ref{thm::asymbehaviour}.(2), there exist $\underline{u}_1,\ldots,\underline{u}_s\in\Z^r$ and $g_1,\ldots,g_s\in\mathbb{Z}$ such that
    \begin{equation}\label{v-min-lin}
        v(R/{\bf I}^{\underline{n}}) = \min\{\underline{u}_j\cdot\underline{n}+g_j : 1\le j\le s\} \ \mbox{ for all } \underline{n}\gg \underline{0},
    \end{equation}
    where the $i$th component $u_{ji}$ of the coefficient vector $\underline{u}_j$ lies in $\{d_{i,1},\dots,d_{i,a_i}\}$ for $1\le i\le r$. In particular, $\underline{u}_j\ge \underline{\delta}$ for $1\le j\le s$.
    Hence, since $v(R/{\bf I}^{\underline{n}+\underline{1}})\le v({\bf I}^{\underline{n}}/{\bf I}^{\underline{n}+\underline{1}})$ for all $\underline{n}\in\N^r$ (cf.~\cite[Prop.~2.5.(2)]{FG}), there exist $g,h\in\Z$ such that
    \begin{equation}\label{eqn:squeeze:v-num}
        \underline{\delta} \cdot \underline{n} + g \le v(R/{\bf I}^{\underline{n}}) \le \underline{\delta} \cdot \underline{n} + h
        \quad \text{for all } \underline{n}\gg \underline{0}.
    \end{equation}
    Following the arguments as shown in the proof of Theorem~\ref{thm:linearity}, one obtains that $v(R/{\bf I}^{\underline{n}})$ is eventually linear with the leading coefficients given by $\underline{\delta}$.
\end{proof}

\section{Examples}\label{sec:Examples}

Here we present several examples that complement our main results. Computations using Macaulay2 \cite{M2} were helpful in constructing some of these examples. For the reader’s convenience, we describe the Macaulay2 commands used.

\begin{para}
    Using Lemma~\ref{lem:$v$-num-indeg}, we only need to compute generators of (usually complicated) modules. This can be done using the Macaulay2 command \texttt{mingens}, which also orders the generators in increasing degree. Other commands we use include \texttt{ass}, which computes the associated prime ideals of modules over a polynomial ring or over quotients of a polynomial ring by a homogeneous ideal, and \texttt{saturate}, which computes modules of the form $I^mJ^n :_R \mathfrak{m}^\infty$. For modules over other rings, for example in \ref{example-5}, the command \texttt{primaryDecomposition} is required, since the command \texttt{ass} does not return any result. By applying \texttt{radical} to the output of \texttt{primaryDecomposition}, one can compute the associated primes.
\end{para}

In the following example, over a polynomial ring on an integral domain, the linearity of $v(R/{\bf I}^{\underline{n}})$ contrasts sharply with the non-linear behaviour of $\reg(R/{\bf I}^{\underline{n}})$.

\begin{example}\label{example-1}
    Let $R=K[x,y]$ be a standard graded polynomial ring in two variables $x$ and $y$ over a field $K$. Set $I:=(x,y^2)$, $J:=(x^2,y)$, and $\fm:=(x,y)$. Then, for all $m,n\in\N$ with $m+n\ge 1$, the following hold.
    \begin{enumerate}[\rm (1)]
        \item $\Ass_R(R/I^m J^n) = \{\fm\}$ and $v(R/I^m J^n) = v_{\fm}(R/I^m J^n) = m+n$.
        \item \cite[Ex.~3.1]{BC17} $\reg(R/I^m J^n) = \max\{m+2n-1, 2m+n-1\}$ .
    \end{enumerate}
\end{example}

\begin{proof}
    Fix $m,n\in\N$ not both zero. Since $x,y\in \sqrt{I^{m}J^{n}}$, $\Ass_R(R/I^m J^n) = \{\fm\}$. Note that the ideal $I^m J^n = (x,y^2)^m (x^2,y)^n$ is given by
    \begin{align*}
        (x^m&,x^{m-1}y^2,x^{m-2}y^{4},\ldots,xy^{2m-2},y^{2m}) (x^{2n},x^{2n-2}y,x^{2n-4}y^2,\ldots,x^2y^{n-1},y^n) \\
        =(&x^{m+2n},x^{m+2n-2}y,x^{m+2n-4}y^2,\ldots,x^{m+2}y^{n-1},x^my^n,\\
        &x^{m+2n-1}y^2,x^{m+2n-3}y^3,x^{m+2n-5}y^4,\ldots,x^{m+1}y^{n+1},x^{m-1}y^{n+2},\ldots,\\
        &x^{2n}y^{2m},x^{2n-2}y^{2m},x^{2n-4}y^{2m+2},\ldots,x^2y^{2m+n-1},y^{2m+n}).
    \end{align*}
    
    Clearly, $\fm = (I^m J^n :_R x^{m-1}y^{n+1}) = (I^m J^n :_R x^{m+1}y^{n-1})$, and $m+n$ is the least possible degree of a homogeneous element of $R/I^m J^n$ whose annihilator is $\fm$. So the assertion in (1) follows. For the equality in (2), note that $\reg(R/I^m J^n) = \reg(I^m J^n)-1 = \max\{m+2n-1, 2m+n-1\}$ by \cite[Ex.~3.1]{BC17}.
\end{proof}

In the next example, none of $v(R/I^m J^n)$ and $v(I^m J^n/I^{m+1} J^{n+1})$ are eventually linear in $(m,n)$. Here, we use the notation $\eend(M) := \sup\{n : M_n\neq 0\}$, where $M$ is a nonzero graded $R$-module.

\begin{example}\label{example-2}
    Let $K[X,Y]$ be a standard graded polynomial ring in two variables $X$ and $Y$ over a field $K$. Set $R:=K[X,Y]/(XY)$, and denote the images of $X$ and $Y$ in $R$ as $x$ and $y$ respectively. Then $R=K[x,y]$. Set $I:=(x,y^2)$, $J:=(x^2,y)$, and $\fm:=(x,y)$. Then, $(0:_RI)=0$ and $(0:_RJ)=0$. Moreover,
    \begin{enumerate}[\rm (1)]
        \item $\Ass_R(R/I^m J^n) = \{\fm\} = \Ass_R(I^{m-1} J^{n-1}/I^m J^n)$ whenever $m, n\ge 1$.
        \item $v(R/I^m J^n) = \min\{ m+2n-1,2m+n-1 \}$ for all $m,n\in\N$ with $m+n\ge 2$.
        \item $\reg(R/I^m J^n) = \max\{m+2n-1,2m+n-1\}$ for all $m,n\in\N$ with $m+n\ge 1$.
        \item $v(I^{m-1} J^{n-1}/I^m J^n) = v(R/I^m J^n)$ for all $m, n\ge 1$.
    \end{enumerate}
\end{example}

\begin{proof}
    Fix $m,n\in\N$ not both zero. It follows
    $$I^m J^n = (x,y^2)^m (x^2,y)^n = (x^m,y^{2m})(x^{2n},y^n) = (x^{m+2n},y^{2m+n}).$$    
    Then, $\Ass_R(R/I^m J^n) = \{\fm\}$. Since $IJ$ annihilates $I^{m-1} J^{n-1}/I^m J^n$, every associated prime of this module will contain $IJ$, and hence must be the same as $\fm$. So (1) follows. Considering the gradation of $R/I^m J^n$, since $R/I^m J^n$ has finite length, $\reg(R/I^m J^n) = \eend(R/I^m J^n) = \max\{m+2n-1,2m+n-1\}$ whenever $m+n\ge 1$. It shows (3). When $m+n\ge 2$, one has that $\fm = (I^m J^n :_R x^{m+2n-1}) = (I^m J^n :_R y^{2m+n-1})$. Moreover, there is no other homogeneous element $f\in R$ of degree different from $m+2n-1$ and $2m+n-1$ such that $\fm=(I^m J^n :_R f)$. Thus, (2) follows. For (4), observe that the images of $x^{m+2n-1}$ and $y^{2m+n-1}$ in $I^{m-1} J^{n-1}/I^m J^n$ are nonzero elements, where $m, n\ge 1$.
\end{proof}

In the next example, both $v_{\fm}(M/I^{m}J^{n}M)$ and $v(M/I^{m}J^{n}M)$ are asymptotically not linear in $(m,n)$. Moreover, in this example, all four functions induced by the local and global $v$-numbers are asymptotically distinct functions.

\begin{example}\label{example-3}
    Let $R=K[X,Y,Z]$ be a standard graded polynomial ring in three variables over a field $K$. Consider the module $M:=R/(XY)$, and the ideals $I:=(X,Z^2)$, $J:=(Y,Z^3)$, $\fp:=(X,Z)$, $\fq:=(Y,Z)$ and $\fm:=(X,Y,Z)$. Then, $(0:_MI)=0$ and $(0:_MJ)=0$. Moreover, for all $m,n\ge 1$, the following hold.
    \begin{enumerate}[\rm (1)]
        \item $\Ass_R(M/I^m J^nM) = \{\fp,\fq,\fm\}$.
        \item $v_{\fm}(M/I^{m}J^{n}M)=\min\{2m+n+1,m+3n\}$.
        \item $v_{\fp}(M/I^{m}J^{n}M)=2m+n-1$ and $v_{\fq}(M/I^{m}J^{n}M)=m+3n-1$.
        \item $v(M/I^{m}J^{n}M)=\min\{2m+n-1,m+3n-1\}$.
    \end{enumerate}
\end{example}

\begin{proof}
    Fix $m,n\ge 1$. The ideal $I^m J^n+(XY)$ is generated by
    \begin{align*}
        XY,&X^mZ^{3n},X^{m-1}Z^{3n+2},X^{m-2}Z^{3n+4},\dots,XZ^{2m+3n-2},\\
        &Y^nZ^{2m},Y^{n-1}Z^{2m+3},Y^{n-2}Z^{2m+6},\dots,YZ^{2m+3n-3},Z^{2m+3n}.
    \end{align*}    
    Since $M/I^mJ^nM\cong R/(I^mJ^n+(XY))$, considering the primary decomposition of the monomial ideal $I^m J^n+(XY)$, one obtains (1).
    
    Denote the images of $X$, $Y$ and $Z$ in $M$ as $x$, $y$ and $z$ respectively. Then $xy=0$. Moreover, the $R$-module $I^mJ^nM$ is generated by
    \begin{align*}
        &x^mz^{3n},x^{m-1}z^{3n+2},x^{m-2}z^{3n+4},\dots,xz^{2m+3n-2},\\
        &y^nz^{2m},y^{n-1}z^{2m+3},y^{n-2}z^{2m+6},\dots,yz^{2m+3n-3},z^{2m+3n}.
    \end{align*}    
    Therefore, $\fm=(I^mJ^n M:_R y^{n-1}z^{2m+2})=(I^mJ^nM:_R x^{m-1}z^{3n+1})$. On the other hand, there are no elements of smaller degree in $M/I^mJ^nM$ whose annihilator is $\fm$. Therefore, $v_{\fm}(M/I^{m}J^{n}M)=\min\{2m+n+1,m+3n\}$, which shows (2).

    In the same manner, one sees that 
    \begin{itemize}
        \item an element which realizes $v_\fp(M/I^{m}J^{n}M)$ is $y^nz^{2m-1}$;
        \item an element which realizes $v_\fq(M/I^{m}J^{n}M)$ is $x^mz^{3n-1}$.
    \end{itemize}
    This implies (3) and (4), which completes the proof.
\end{proof}

\begin{remark}
    In Example~\ref{example-3}, both $v_\fp(M/I^mJ^nM)$ and $v_\fq(M/I^mJ^nM)$ eventually become linear, but $v(M/I^mJ^nM)=\min\{v_\fp(M/I^mJ^nM),v_\fq(M/I^mJ^nM)\}$ is not eventually linear. Moreover, asymptotically, $v(M/I^{m}J^{n}M)$ and $v_{\fm}(M/I^{m}J^{n}M)$ are two different functions.
\end{remark}

The following example ensures that, despite Theorem~\ref{thm:lin-v-function}, one cannot expect that every local $v$-number for products and powers of several ideals eventually becomes linear even over a polynomial ring over a field.

\begin{example}\label{example-4}
    Let $R=K[x,y,z]$ be a standard graded polynomial ring over a field $K$. Consider the ideals $I=(x^2,yz^2)$ and $J=(y^2,xz^2)$. Set $\fp:=(x,y)$, $\fq :=(x,z)$, $\fr :=(y,z)$, and $\fm:=(x,y,z)$. Then, for every $m,n\ge 1$, the following hold.
    \begin{enumerate}[\rm (1)]
        \item $\Ass_R(R/I^m J^n)=\{\fp,\fq,\fr,\fm\}$.
        \item $v_{\fm}(R/I^{m}J^{n})=2m+2n+2$.
        \item $v_{\fp}(R/I^{m}J^{n})=\min\{3m+2n+1,2m+3n+1\}$.
        \item $v(R/I^{m}J^{n})=v_{\fq}(R/I^{m}J^{n})=v_{\fr}(R/I^{m}J^{n})=2m+2n+1$.
        \item $\Ass_R(I^{m-1} J^{n-1}/I^m J^n)=\{\fp,\fq,\fr,\fm\}$, and the (local) $v$-numbers of the modules $I^{m-1}J^{n-1}/I^mJ^n$ and $R/I^mJ^n$ coincide.
    \end{enumerate}
\end{example}

\begin{proof}
    Fix $m,n\ge 1$. The ideals $I^m$ and $J^n$ are generated by
    \begin{align*}
        &\{(x^2)^{s}(yz^2)^{m-s} = x^{2s}y^{m-s}z^{2m-2s} : 0\le s\le m\} \mbox{ and}\\
        &\{(y^2)^{t}(xz^2)^{n-t} = x^{n-t}y^{2t}z^{2n-2t} : 0 \le t \le n\} \mbox{ respectively}.
    \end{align*}
    Hence $I^mJ^n = \big(x^{n+2s-t}y^{m-s+2t}z^{2(m+n-s-t)} : 0\le s\le m, 0 \le t \le n\big)$. The statement in (1) follows from the primary decomposition of this monomial ideal. For better understanding, the reader may consider the case $m,n=1$.
    
    In view of Lemma~\ref{lem:$v$-num-indeg}, the local $v$-numbers are given by
    \begin{align*}
        v_\fm(R/I^mJ^n)&=\indeg\big((I^mJ^n:_R\fm)/I^mJ^n\big) \ \mbox{ and}\\
        v_\fa(R/I^mJ^n)&=\indeg\left(\frac{(I^mJ^n:_R\fa)}{(I^mJ^n:_R\fa)\cap (I^mJ^n:_R\fm^\infty)}\right) \mbox{ for } \fa \in \{\fp,\fq,\fr\},
    \end{align*}
    where $(I^mJ^n:_R\fm^\infty) = \bigcup_{\ell\ge 1} (I^mJ^n:_R\fm^\ell)$. Note that
    \begin{align*}
        (I^mJ^n:_R\fp)\cap (I^mJ^n:_R\fm^\infty) = (I^mJ^n:_R\fp)\cap (I^mJ^n:_Rz^\infty), \\
        (I^mJ^n:_R\fq)\cap (I^mJ^n:_R\fm^\infty) = (I^mJ^n:_R\fq)\cap (I^mJ^n:_Ry^\infty), \\
        (I^mJ^n:_R\fr)\cap (I^mJ^n:_R\fm^\infty) = (I^mJ^n:_R\fr)\cap (I^mJ^n:_Rx^\infty).
    \end{align*}
    The assertions (2), (3) and (4) can be obtained from the following observations:
    \begin{align}
        x^{2m-1}y^{2n}z^3, x^{2m}y^{2n-1}z^3 & \in (I^mJ^n:_R\fm) \smallsetminus I^mJ^n, \label{colon-m}\\
        y^{m+2n-1}z^{2m+2}, x^{2m+n-1}z^{2n+2} & \in (I^mJ^n:_R\fp) \smallsetminus (I^mJ^n:_R z^\infty), \label{colon-p}\\
        x^{2m-1}y^{2n+1}z & \in (I^mJ^n:_R\fq) \smallsetminus (I^mJ^n:_Ry^\infty), \label{colon-q}\\
        \mbox{and} \quad x^{2m+1}y^{2n-1}z & \in (I^mJ^n:_R\fr) \smallsetminus (I^mJ^n:_Rx^\infty).\label{colon-r}
    \end{align}    
    These are the monomials of the minimum possible degree contained in the right-hand sides above, and therefore they compute the respective local $v$-numbers.

    The monomials listed in \eqref{colon-m}, \eqref{colon-p}, \eqref{colon-q} and \eqref{colon-r} all lie in $I^{m-1}J^{n-1}$. The annihilator ideals of these monomials in $I^{m-1}J^{n-1}/I^mJ^n$ provide the associated prime ideals $\{\fp,\fq,\fr,\fm\}$ of the module $I^{m-1}J^{n-1}/I^mJ^n$. Due to the minimality of the degrees, these monomials also compute the respective local $v$-numbers of the same module, which coincide with that of $R/I^mJ^n$. Thus, assertion (5) follows.
\end{proof}

\begin{remark}
    The assumption ${\bf I}^{\underline{s}} M\subseteq N$ for some $\underline{s}\in\N^r$, in Theorem~\ref{thm::asymbehaviour}.(2) and Corollary~\ref{cor:asymp-v-num}, seems to be necessary because in the critical case ${\bf I}\subseteq \sqrt{\ann_R(N)}$, the function $v_\fp(M/{\bf I}^{\underline{n}} N)$ is eventually constant for each $\fp\in\mathcal{A}_N^M({\bf I})$. However, the authors are unaware of any example in which ${\bf I}^{\underline{n}} M\not\subseteq N$ for every $\underline{n}\in\N^r$, ${\bf I}\not\subseteq \sqrt{\ann_R(N)}$, and $v_\fp(M/{\bf I}^{\underline{n}} N)$ eventually is not the minimum of linear functions. Some attempts of explicit computations have been made in the case $r=1$. As shown in \cite{FG}, the same hypotheses lead to a linear behaviour. Due to these attempts, we feel that the condition ${\bf I}^{\underline{s}} M\subseteq N$, which is merely technical but still true in many cases, could be relaxed.
\end{remark}

The following example shows that the condition ``$R_0$ is an integral domain" in Theorem~\ref{thm:lin-v-function} is crucial.

\begin{example}\label{example-5}
    Set $R_0:=K[z,w]/(zw)$, where $z,w$ are indeterminates over a field $K$. Let $R=R_0[x,y]$ be a polynomial ring over $R_0$, and consider the homogeneous ideals $I=(x^2,wx)$ and $J=(y^2,zy)$. Then, for any $m,n\ge 1$, one has
    \[
    \indeg(I^mJ^n/I^{m+1}J^{n+1})=\min\{2m+n,m+2n\}.
    \]
\end{example}

\begin{proof}
     Fix $m,n\ge 1$. The ideals $I^m$ and $J^n$ are generated by
    \begin{align*}
        &\{(x^2)^{s}(wx)^{m-s} = x^{s+m}w^{m-s} : 0\le s\le m\} \mbox{ and}\\
        &\{(y^2)^{t}(zy)^{n-t} = y^{t+n}z^{n-t} : 0 \le t \le n\} \mbox{ respectively}.
    \end{align*}
    Since $zw=0$ in $R_0$, it follows that the product $I^mJ^n$ is generated by
    \[
    \big(x^{2m}y^{t+n}z^{n-t}, x^{s+m}y^{2n}w^{m-s}: 0\le s\le m, 0 \le t \le n\big).
    \]
    We note that the elements in $R_0$ have zero degree. Moreover, we can see that $x^{2m}y^n$ and $x^my^{2n}$ do not lie in $I^{m+1}J^{n+1}$. Therefore, the result follows.
\end{proof}

\subsection*{Acknowledgments}
The authors thank Aldo Conca for his advice on this paper and for encouraging this collaboration, and Maria Evelina Rossi for the helpful discussions on this problem. The first named author is partially supported by INdAM-GNSAGA, MIUR Excellence Department Project CUP\, D33C23001110001, and by PRIN 2022 Project 2022K48YYP. The authors also thank the anonymous referee for carefully reading the article and for providing many helpful comments and suggestions that have improved its presentation.

\section*{\bf Declarations}
{\bf Ethical approval:} The work in this article is original and is not submitted anywhere else for possible publication.

{\bf Competing interest:} None of the authors have competing interests.

{\bf Authors' contributions:} Both the authors have equal contributions on each of the following: Investigation, Methodology, Writing the article and Review.

{\bf Funding:} Fiorindo is partially supported by INdAM-GNSAGA, MIUR Excellence Department Project CUP\, D33C23001110001, and by PRIN 2022 Project 2022K48YYP.

{\bf Availability of data and materials:} This article has no associated data.


\begin{thebibliography}{AAAA}


\bibitem{MR1287608} W.W.~Adams and P.~Loustaunau, {\it An introduction to {G}r\"obner bases}, Graduate Studies in Mathematics \textbf{3}. American Mathematical Society, 1994.

\bibitem{BM} P.~Biswas and M.~Mandal, {\it A study of 
v-number for some monomial ideals}, Collect. Math. (2024). 

\bibitem{brod} M.~Brodmann, {\it Asymptotic stability of ${\rm Ass}(M/I\sp{n}M)$}, Proc. Amer. Math. Soc. \textbf{74} (1979), 16--18.

\bibitem{BC17} W.~Bruns and A.~Conca, {\it A remark on regularity of powers and products of ideals}, J. Pure Appl. Algebra {\bf 221} (2017), 2861--2868.

\bibitem{BCRV} W.~Bruns, A.~Conca, C.~Raicu and M.~Varbaro, {\it Determinants, Gröbner bases and cohomology}, Springer Monogr. Math. (2022).

\bibitem{Conca23} A.~Conca, {\it A note on the $v$-invariant}, Proc. Amer. Math. Soc. {\bf 152(6)} (2023), 2349--2351.

\bibitem{CSTPV} S. M.~Cooper, A.~Seceleanu, S.O.~Toh\v{a}neanu, M.~Pinto Vaz and R.H.~Villarreal, {\it Generalized minimum distance functions and algebraic invariants of Geramita ideals}, Adv.~Appl. Math. {\bf 112} (2020), 101940, 34 pp.

\bibitem{FS} A.~Ficarra and E.~Sgroi, {\it Asymptotic behaviour of the v-number of homogeneous ideals}, \href{https://arxiv.org/abs/2306.14243}{arXiv:2306.14243}.

\bibitem{FG} L.~Fiorindo and D.~Ghosh, {\it On the asymptotic behaviour of the Vasconcelos invariant for graded modules}, Nagoya Math. J. (2025), 1--15.

\bibitem{Gh16} D.~Ghosh, {\it Asymptotic linear bounds of Castelnuovo--Mumford regularity in multigraded modules}, J.~Algebra {\bf 445} (2016), 103--114.

\bibitem{GP25} D.~Ghosh and S.~Pramanik, {\it Asymptotic v-numbers of graded (co)homology modules involving powers of an ideal}, J.~Algebra {\bf 671} (2025), 61--74.
    
\bibitem{M2} D. R.~Grayson and M. E.~Stillman, {\it Macaulay2, a software system for research in algebraic geometry}, \href{http://www2.macaulay2.com}{http://www2.macaulay2.com}.
    


\bibitem{Hay06} F.~Hayasaka, {\it Asymptotic stability of primes associated to homogeneous components of multigraded modules}, J.~Algebra \textbf{306} (2006), 535--543.

\bibitem{JS} D.~Jaramillo-Velez and L.~Seccia, {\it Connected domination in graphs and v-numbers of binomial edge ideals}, Collect. Math. {\bf 75} (2024), no.~3, 771--793.

\bibitem{JV} D.~Jaramillo and R.~Villarreal, {\it The v-number of edge ideals}, Journal Of Combinatorial Theory, Series A \textbf{177}, Paper No. 105310 (2021).


\bibitem{KSharp} A.~Kingsbury and R.~Sharp, {\it Asymptotic behaviour of certain sets of prime ideals}, Proc. Amer. Math. Soc. \textbf{124} (1996), 1703--1711.

\bibitem{KR} M.~Kreuzer and L.~Robbiano, {\it Computational commutative algebra. 1}, Springer, Berlin, 2000; MR1790326

\bibitem{Sa} K.~Saha, {\it The v-Number and Castelnuovo-Mumford Regularity of Cover Ideals of Graphs}, Int. Math. Res. Not. (2023), 9010--9019.

\bibitem{NS} K.~Saha and N.~Kotal, {\it On the v-number of Gorenstein ideals and Frobenius powers}, Bull. Malays. Math. Sci. Soc. {\bf 47} (2024), no.~6, Paper No. 167.


\bibitem{SS} K.~Saha and I.~Sengupta, {\it The v-number of monomial ideals}, J. Algebraic Combin. \textbf{56} (2022), 903--927.

\bibitem{West04} E.~West, {\it Primes associated to multigraded modules}, J.~Algebra {\bf 271} (2004), 427--453.

\end{thebibliography}
\end{document}